\documentclass[11pt]{article}
\usepackage{amssymb}
\usepackage{amsmath}
\usepackage{latexsym}

\input{prepictex}
\input{pictex}
\input{postpictex}
\begin{document}
\bibliographystyle{plain}
\baselineskip=18pt
\newtheorem{lemma}{Lemma}[section]
\newtheorem{theorem}[lemma]{Theorem}
\newtheorem{prop}[lemma]{Proposition}
\newtheorem{cor}[lemma]{Corollary}
\newtheorem{definition}[lemma]{Definition}
\newtheorem{definitions}[lemma]{Definitions}
\newtheorem{remark}[lemma]{Remark}
\newtheorem{conjecture}[lemma]{Conjecture}
\newcommand{\TT}[1]{\widetilde{T}\mbox{}_{#1}}
\newcommand{\bb}{{\bf B}'({\bf d})\cup {\bf B}''({\bf d})}
\newcommand{\vanish}{{\bf B}\setminus \Psi(\bb)}
\newcommand{\word}[1]{s_{i_1}s_{i_2}\cdots s_{i_{#1}}}
\newcommand{\dd}{\mbox{\scriptsize\bf d}}
\newcommand{\E}[1]{E_{i_1}E_{i_2}\cdotsE_{i_{#1}}}
\newcommand{\F}[1]{F_{i_1}F_{i_2}\cdotsF_{i_{#1}}}
\newcommand{\EA}[1]{E_{i_1}^{(r_1)}E_{i_2}^{(r_2)}\cdots
           E_{i_{#1}}^{(r_{#1})}}
\newcommand{\FA}[1]{F_{i_1}^{(r_1)}F_{i_2}^{(r_2)}\cdots
           F_{i_{#1}}^{(r_{#1})}}
\newcommand{\UA}{U_{\cal A}^-}
\newcommand{\dimr}{\left( \begin{array}{c} n+1 \\ r \end{array} \right)}
\newcommand{\TE}{\widetilde{E}}
\newcommand{\TF}{\widetilde{F}}
\newcommand{\C}{{\mathbb{C}}}
\newcommand{\N}{{\mathbb{N}}}
\newcommand{\Q}{{\mathbb{Q}}}
\newcommand{\R}{{\mathbb{R}}}
\newcommand{\Z}{{\mathbb{Z}}}

\newcommand{\B}{{\mathbf{B}}}

\newcommand{\CA}{{\mathcal{A}}}
\newcommand{\CL}{{\mathcal{L}}}

\newcommand{\0}{-\,}

\newcommand{\cdi}{\mbox{CD$(\mathbf{i})$}}
\newcommand{\cdii}{\mbox{CD$(\mathbf{i}')$}}

\begin{center}
\baselineskip=13pt
{\LARGE The Lusztig cones of a quantized enveloping algebra of type $A$}

\vspace*{10mm}
{\Large Robert Marsh}

{\em Department of Mathematics and Computer Science, University of Leicester,
University Road, Leicester LE1 7RH, England} \\
{E-mail: R.Marsh@mcs.le.ac.uk} \\
{Fax: +44 (0) 116 252 3915}

\baselineskip=21pt

{\bf Abstract} \\
\parbox[t]{5.15in}{
We show that for each reduced expression for the longest word in the Weyl
group of type $A_n$, the corresponding cone arising in Lusztig's description
of the canonical basis in terms of tight monomials is simplicial,
and construct explicit spanning vectors.

} \\ \ \\
\end{center}

\section{Introduction}

Let $U=U_q(sl_{n+1}(\mathbb{C}))$ be the quantum group associated to the
simple Lie algebra $\mathfrak{g}=sl_{n+1}(\mathbb{C})$ of type $A_n$.
The negative part $U^-$ of $U$ has
a canonical basis with favourable properties
(see Kashiwara~\cite{kash2} and Lusztig~\cite[\S14.4.6]{lusztig6}).
For example, via action on highest weight vectors it gives rise to bases
for all the finite-dimensional irreducible highest weight $U$-modules.

The Lusztig cones first appeared in a paper~\cite{lusztig7} of Lusztig,
where he proved that, if $\mathbf{i}=(i_1,i_2,\ldots ,i_k)$ is a reduced
expression for the longest word $w_0$ in the Weyl group of $\mathfrak{g}$
(i.e. $s_{i_1}s_{i_2}\cdots s_{i_k}$ is reduced), and
$\mathbf{a}=(a_1,a_2,\ldots ,a_k)$ satisfies certain linear inequalities,
then the monomial
$$F_{i_1}^{(a_1)}F_{i_2}^{(a_2)}\cdots F_{i_k}^{(a_k)}$$
lies in the canonical basis of $U^-$, provided $\mathfrak{g}$ is of type
$A_1$, $A_2$ or $A_3$. We call the set of points associated to
$\mathbf{i}$ satisfying these linear equalities the {\em Lusztig cone}
associated to $\mathbf{i}$.
This result was extended to type $A_4$ in~\cite{me7}, but
has been shown to fail for larger $n$ by Reineke~\cite{reineke2} and
Xi~\cite{xi4}. However, in the crystal limit of Kashiwara, this cone plays an
important role. This role will be explored in~\cite{me9}, and the joint
paper~\cite{me10} with Roger Carter, where a link is made with the regions
of linearity of Lusztig's reparametrization functions. Indeed, it may be
that the Lusztig cones are regions of linearity for
reparametrization functions arising from a change of reduced expression for
$w_0$ in the string parametrization of the canonical basis.

The Lusztig cones have also arisen in other contexts. For example, let
$G=SL_n(\mathbb{C})$, and let $x_1,x_2,\ldots ,x_n$ be a choice of
one-parameter simple root subgroups of $G$. Zelevinsky~\cite{zelevinsky2}
shows that a generic element $x\in L^{e,w_0}$, a reduced real double Bruhat
cell of $G$, has form $x=x_{i_1}(t_1)x_{i_2}(t_2)\cdots x_{i_k}(t_k)$, where
$\mathbf{i}$ is any fixed reduced expression for $w_0$. The parameters
$t_1,t_2,\ldots ,t_k$ can be regarded as functions on $L^{e,w_0}$, and
Zelevinsky shows that the Laurent monomial $t_1^{a_1}t_2^{a_2}\cdots
t_k^{a_k}$ is regular if and only if $\mathbf{a}=(a_1,a_2,\ldots ,a_k)$
satisfies the defining
inequalities of the Lusztig cone corresponding to $\mathbf{i}$.
Lusztig cones also play a role in the paper~\cite{nz1}
investigating the string cones --- cones which arise from the string
parametrization of the canonical basis of $U^-$ (which also depends on a
choice of reduced expression $\mathbf{i}$ for $w_0$).

In this paper, we investigate the Lusztig cones and show that they possess a
beautiful combinatorial structure. We prove that each Lusztig
cone is the set of integral points of a simplicial cone
and give an explicit description of spanning vectors for it, in terms
of the corresponding reduced expression for $w_0$. Provided we label points
in a Lusztig cone in the correct way, in this description there are $n$
spanning vectors common to each Lusztig cone. The other vectors correspond
to the bounded chambers of the `chamber ansatz' for the
corresponding reduced expression (see~\cite[\S\S1.4, 2.3]{bfz1}),
and each vector depends only on the set of braids which pass below the
corresponding chamber.

We give an explicit description of these spanning vectors, for an arbitrary
reduced expression for $w_0$, using a relabelling of the bounded chambers
with `partial quivers' --- certain graphs with a mixture of oriented and
unoriented edges such that the subgraph of oriented edges is connected. These
allow us to give a natural description of the spanning vectors
(Theorem~\ref{maintheorem}),
and also enable us to use the quiver-compatible reduced expressions for $w_0$
as a way of approaching all reduced expressions.

In the case where $\mathbf{i}$ is quiver-compatible, for simply-laced type,
these vectors have now also been studied by Bedard
in~\cite{bedard2}. Bedard describes these vectors
using the Auslander-Reiten quiver of the quiver and homological algebra,
showing they are closely connected to the representation theory of the quiver.

The description of the Lusztig cones given here should be
of use in understanding the canonical basis. It is known that the Lusztig
cone is contained in the corresponding string cone --- this is shown
by Premat~\cite{premat1} (a proof which holds for all symmetrizable Kac-Moody
Lie algebras) and also in~\cite{me9} for type $A_n$ using the results derived
here. This means that there is a subset of the canonical basis corresponding
to each Lusztig cone and therefore a canonical basis element corresponding to
each spanning vector described here. In types $A_1$, $A_2$ and $A_3$, the
corresponding dual canonical basis element is known to be primitive --- that
is, it cannot be factorized as a product of dual canonical basis elements
(this can be seen using~\cite{bz1}). These primitive elements, and whether
they quasi-commute, play a key role in current understanding of the dual
canonical basis (see~\cite{bz1},~\cite{caldero1} and~\cite{lz1}),
and also in the structure of the module categories of affine Hecke algebras
(see~\cite{lnt1}).

\section{Lusztig cones are simplicial} \label{spanningvectors}

For positive integers $p<q$ we denote by $[p,q]$ the set
$\{p,p+1,\ldots ,q\}$, and for a rational number $x$ we denote by
$\lceil x \rceil$ the smallest element of $\{y\in \mathbb{Z}\,:\,x\leq y\}$.
Let $\mathfrak{g}$ be the simple Lie algebra $sl_{n+1}(\mathbb{C})$,
with root system $\Phi$, and simple roots
$\alpha_1,\alpha_2,\ldots ,\alpha_n$. We use the following numbering of the
Dynkin diagram (and its edges):


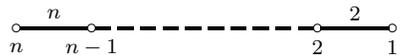
\beginpicture

\setcoordinatesystem units <1cm,1cm>             
\setplotarea x from -1 to 12, y from 1.5 to 3       

\scriptsize{

\multiput {$\circ$} at 3   2 *1 1 0 /      %
\multiput {$\circ$} at 7   2 *1 1 0 /      

\linethickness=1pt           

\putrule from 3.05 2 to 3.95 2  %
\putrule from 7.05 2 to 7.95 2  

\setdashes <2mm,1mm>          %
\putrule from 4.05 2 to 6.95 2  

\put {$n$}   [c] at 3 1.75
\put {$n-1$}   [c] at 4 1.75
\put {$2$} [c] at 7 1.75
\put {$1$} [c] at 8 1.75

\put{$n$}[c] at 3.5 2.2
\put{$2$}[c] at 7.5 2.2

}

\endpicture


Let $W$ be the Weyl group of $\mathfrak{g}$. Recall that $W$ can be presented
with generators $s_1,s_2,\ldots ,s_n$, subject to relations:
\begin{eqnarray}
s_i^2 & = & e, \\
s_is_j & = & s_js_i \mbox{\ \ \ if\ }|i-j|>1, \label{rel1} \\
s_is_js_i & = & s_js_is_j \mbox{\ \ \ if\ }|i-j|=1. \label{rel2}
\end{eqnarray}
We recall also that relations of type (\ref{rel1}) are called
{\em short braid relations},
and relations of type (\ref{rel2}) are called {\em long braid relations}. Two
reduced expressions $w_0$ are said to be {\em commutation equivalent} if there
exists a sequence of short braid relations taking the first to the second;
the equivalence classes of this equivalence relation are called {\em
commutation classes}.
Let $w_0$ be the longest element in $W$.
We shall identify a reduced expression for $s_{i_1}s_{i_2}\cdots s_{i_k}$
for $w_0$ with the corresponding $k$-tuple ${\bf i}=
(i_1,i_2,\ldots ,i_k)$.
The {\em Lusztig cone}, $C_{\bf i}$, corresponding
to $\bf i$ is defined to be the set of points ${\bf a}\in \mathbb{N}^k$
satisfying the following condition: \\
For every pair $i_s,i_{s'}$, with $s,s'\in
[1,k]$, $i_s=i_{s'}=i$ and such that $i_p\not=i$ whenever
$s<p<s'$, we have
\begin{equation}\label{ineqs}
(\sum_p a_p)-a_s - a_{s'}\geq 0,
\end{equation}
where the sum is over all $p$ with $s<p<s'$ such that $i_p$ is joined to $i$
in the Dynkin diagram.

In this section, we show that the Lusztig cone corresponding to any reduced
expression for the longest word $w_0$ is simplicial.
We can regard $C_{\bf i}$ as the subset of $\mathbb{Z}^k$ defined by
the $k-n$ inequalities~(\ref{ineqs}) above, together with the $k$ inequalities
$a_j\geq 0$ for $j=1,2,\ldots ,k$. We will show that only $n$ of the
inequalities $a_j\geq 0$ are necessary to define $C_{\bf i}$.

The reduced expression $\bf i$ defines an ordering on the set $\Phi^+$ of
positive roots of the root system associated to $W$.
If we write $\alpha^j=s_{i_1}s_{i_2}\cdots s_{i_{j-1}}(\alpha_{i_j})$ for
$j=1,2,\ldots ,k$, then $\Phi^+=\{\alpha^1,\alpha^2,\ldots \alpha^k\}$
(and there are no repetitions in this list). For ${\bf a}=(a_1,a_2,\ldots
,a_k)\in \mathbb{Z}^k$,
write $a_{\alpha^j}=a_j$. We shall use this relabelling
of $\bf a$ throughout.

\begin{figure}
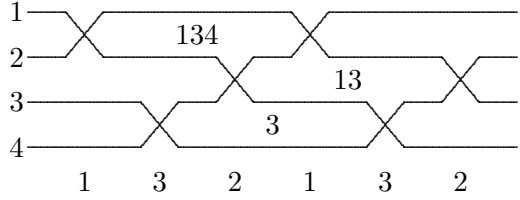

\beginpicture

\setcoordinatesystem units <0.5cm,0.3cm>             
\setplotarea x from -10 to 14, y from -3 to 8       
\linethickness=0.5pt           

\put{$1$}[c] at -0.3 6 %
\put{$2$}[c] at -0.3 4 %
\put{$3$}[c] at -0.3 2 %
\put{$4$}[c] at -0.3 0 %

\put{$1$}[c] at 1.5 -1.5 %
\put{$3$}[c] at 3.5 -1.5 %
\put{$2$}[c] at 5.5 -1.5 %
\put{$1$}[c] at 7.5 -1.5 %
\put{$3$}[c] at 9.5 -1.5 %
\put{$2$}[c] at 11.5 -1.5 %

\setlinear \plot 0 0  3 0 / %
\setlinear \plot 3 0  4 2 / %
\setlinear \plot 4 2  5 2 / %
\setlinear \plot 5 2  6 4 / %
\setlinear \plot 6 4  7 4 / %
\setlinear \plot 7 4  8 6 / %
\setlinear \plot 8 6  13 6 / %

\setlinear \plot 0 2  3 2 / %
\setlinear \plot 3 2  4 0 / %
\setlinear \plot 4 0  9 0 / %
\setlinear \plot 9 0  10 2 / %
\setlinear \plot 10 2  11 2 / %
\setlinear \plot 11 2  12 4 / %
\setlinear \plot 12 4  13 4 / %

\setlinear \plot 0 4  1 4 / %
\setlinear \plot 1 4  2 6 / %
\setlinear \plot 2 6  7 6 / %
\setlinear \plot 7 6  8 4 / %
\setlinear \plot 8 4  11 4 / %
\setlinear \plot 11 4  12 2 / %
\setlinear \plot 12 2  13 2 / %

\setlinear \plot 0 6  1 6 / %
\setlinear \plot 1 6  2 4 / %
\setlinear \plot 2 4  5 4 / %
\setlinear \plot 5 4  6 2 / %
\setlinear \plot 6 2  9 2 / %
\setlinear \plot 9 2  10 0 / %
\setlinear \plot 10 0  13 0 / %

\put{$134$}[c] at 4.5 5 %
\put{$13$}[c] at 8.5 3 %
\put{$3$}[c] at 6.5 1 %
\endpicture
\caption{The Chamber Ansatz of $(1,3,2,1,3,2)$.\label{chamberdiagram}}
\end{figure}

We shall need the chamber ansatz for $\mathbf{i}$
defined in~\cite[\S\S1.4, 2.3]{bfz1}.
We take $n+1$ strings, numbered from top to bottom, and write 
$\mathbf{i}$ from left to right along the bottom of the diagram.
Above a letter $i_j$ in $\mathbf{i}$,
the $i_j$th and $(i_j+1)$st strings from the top above $i_j$ cross.
Thus, for example, in the case $n=3$ with $\mathbf{i}=(1,3,2,1,3,2)$, the
chamber ansatz is shown in Figure~\ref{chamberdiagram}.
We denote the chamber ansatz of $\mathbf{i}$ by CD($\mathbf{i}$).
A {\em chamber} will be regarded as a pair $(c,{\mathbf{i}})$, where $c$ is a
component of the complement of CD($\mathbf{i}$). We shall deal only with
bounded chambers, so in the sequel, `chamber' will always mean bounded
chamber. Each chamber $(c,{\mathbf{i}})$ can be labelled with the numbers
of the strings passing below it, denoted $l(c,{\mathbf{i}})$.
Following~\cite{bfz1}, we call such a label a chamber set.
For example, the chamber sets corresponding to the $3$ chambers in
Figure~\ref{chamberdiagram} are $134,3,13$. Note that the set of chamber
sets of $\bf i$ is independent of its commutation class
(so we can talk also of the chamber sets of such a commutation class).
Every non-empty
subset of $[1,n+1]$ can arise as a chamber set for some $\bf i$,
except for subsets of the form $[1,j]$ and $[j,n+1]$ for $1\leq j\leq n+1$;
this is observed in the proof of~\cite[Theorem 2.7.1]{bfz1}.

Since the crossings of strings correspond to letters in $\bf i$, and the 
horizontal height of a crossing is the value of the corresponding 
letter, the chambers $(c,{\bf i})$ correspond precisely to pairs of equal
letters in $\bf i$ which have no further occurrence of that letter
between them. We call such a pair {\em minimal}.
Next, label the crossing point of strings $i$ and $j$
with the pair $(i,j)$. If we identify $(i,j)$ with the positive root
$\alpha_{ij}=\alpha_i+\alpha_{i+1}+\cdots +\alpha_{j-1}$ then the labels
of the crossings will be precisely the positive roots in the ordering
defined (above) by $\bf i$. We will write $a_{ij}$ for $a_{\alpha_{ij}}$.

The key to understanding $C_{\bf i}$ is the following Lemma.

\begin{lemma} \label{addroots}
Suppose ${\bf a}\in C_{\bf i}$, and that $\alpha$, $\beta$ and
$\alpha+\beta$ are all positive roots. Then $a_{\alpha+\beta}\geq
a_{\alpha}+a_{\beta}$.
\end{lemma}

{\bf Proof:} We note that each chamber $(c,{\bf i})$ determines a minimal
pair in $\bf i$ and thus a defining inequality of $C_{\bf i}$. Such an
inequality can be written
$$(\sum a_{\gamma}) - a_{\gamma_1} - a_{\gamma_2}\geq 0$$
where $\gamma_1$ and $\gamma_2$ are the labels of the crossings at the
left and right hand ends of $c$, respectively, and the sum
is over the positive roots $\gamma$ labelling the other crossings
on the boundary of $(c,{\bf i})$.

Suppose that $\alpha=\alpha_{ij}$, and $\beta=\alpha_{jk}$ for some
$1\leq i<j<k\leq n+1$,
so that $\alpha+\beta=\alpha_{ik}$. By inspection, it can be seen that the
sum of the defining inequalities of $C_{\mathbf{i}}$ corresponding to chambers
inside the triangle formed by the pseudolines $i$, $j$ and $k$ gives
precisely the inequality $a_{ik}\geq a_{ij}+a_{jk}$.~$\square$

It follows from Lemma~\ref{addroots} that $C_{\bf i}$ can be described
as the set of points in $\mathbb{Z}^k$ satisfying the
inequalities~(\ref{ineqs}),
together with the inequalities $a_{\alpha}\geq 0$ for $\alpha$ a simple
root. This is a total of $k-n+n=k$ inequalities. We thus have:

\begin{cor}
Let $\mathbf{i}$ be a reduced expression for the longest word $w_0$ in type
$A_n$. Then the corresponding Lusztig cone $C_{\mathbf{i}}$ is simplicial.
\end{cor}

It follows that there is a $k\times k$ matrix $P_{\bf i}\in M_k(\mathbb{Z})$
such that
$$C_{\bf i}=\{{\bf a}\in \mathbb{Z}^k\,:\,P_{\bf i}{\bf a}\geq 0\},$$
where for ${\bf z}\in \mathbb{Z}^k,$ ${\bf z}\geq 0$
means that each entry in $\bf z$ is
nonnegative. To write down $P_{\bf i}$ we fix an ordering on the set
of $k$ inequalities defining $C_{\bf i}$ --- they are indexed by the simple
roots and the chambers of $\mathbf{i}$.
We do this for each $\bf i$. We also have:

\begin{lemma} \label{unimod}
The matrix $P_{\bf i}\in M_k(\mathbb{Z})$ has an inverse in $M_k(\mathbb{N})$.
\end{lemma}

\noindent {\bf Proof:} Let $\alpha=\alpha_{ij}$ be any positive root. By the proof of
Lemma~\ref{addroots}, we know that a sum of rows of $P_{\bf i}$ gives us
a row with corresponding inequality
$a_{ij}-a_{i,i+1}-a_{i+1,i+2}-\cdots -a_{j-1,j}\geq 0$.
But we also have rows of $P_{\bf i}$ corresponding to $a_{i,i+1}\geq 0$,
$a_{i+1,i+2}\geq 0$,\ldots ,$a_{j-1,j}\geq 0$ (as $\alpha_{i,i+1}$, etc.,
are all simple roots). Adding together all of these rows thus gives us
a row with corresponding inequality $a_{ij}\geq 0$, i.e. a row
$(0,0,\ldots ,0,1,0,0,\ldots ,0)$, where the $1$ appears in the position
corresponding to $a_{ij}$ (that is, in the $l$th position, where
$\alpha^l=\alpha_{ij}$ in the ordering defined by $\bf i$).
We can do this for any $l$, so if we write the coefficients of the $l$th
linear combination as the $l$th row of a matrix $Q_{\bf i}$ then
$Q_{\bf i}P_{\bf i}$ is the identity matrix and we are done.~$\square$

Now let $v_1,v_2,\ldots ,v_k$ be the columns of $Q_{\bf i}$. Because
$Q_{\bf i}$ and $P_{\bf i}$ are inverses we have that ${\bf a}\in
\mathbb{Z}^k$ lies in $C_{\bf i}$ if and only if $P_{\bf i}{\bf a}\geq 0$,
if and only if ${\bf a}$ is a nonnegative integer linear combination of
$v_1,v_2,\ldots ,v_k$. It is now clear that $C_{\bf i}$ is the set of
integral points of the simplicial cone in $\mathbb{R}^k$ defined by the same
inequalities as $C_{\bf i}$.
Our next step in understanding $C_{\bf i}$ is to
describe the `spanning vectors' $v_j$, $j=1,2,\ldots ,k$, for each $\bf i$.
It will turn out that these spanning vectors can be neatly described in
terms of `partial quivers', which we define below, together with some
vectors common to every $C_{\bf i}$. Firstly, we investigate
chambers and chamber sets further.

\begin{definitions} \rm \label{chamberequiv} \ \\
It is well-known that the graph with vertices $\chi_n$, the set of reduced
expressions for $w_0$, where ${\bf i},{\bf i'}\in \chi_n$ are linked by an
edge whenever there is a braid relation taking $\bf i$ to $\bf i'$,
is connected.
Note also that, applying a braid relation to a reduced expression $\bf i$
has a corresponding effect on the chamber ansatz for $\bf i$.
If this is a short braid relation $s_is_j=s_js_i$ with $A_{ij}=0$ then
the only effect on the chamber ansatz for $\bf i$ is to reorder two
crossings which do not interfere with each other.
The chambers of $\bf i$ are clearly mapped bijectively onto
the chambers of $\bf i'$ by this operation, and if $(c,{\bf i})$ is a
chamber, with corresponding chamber $(c',{\bf i'})$ (under this mapping),
we set $(c,{\bf i})\sim (c',{\bf i'})$.
If we have a long braid relation $B\,:\,s_is_js_i=s_js_is_j$
with $A_{ij}=A_{ji}=
-1$, then this determines a chamber in
$CD({\bf i})$, with ends
given by the crossings corresponding to the two $s_i$'s in
$s_is_js_i$. If $(c,{\bf i})$ is any other chamber of $\bf i$ then there
is a corresponding pair of letters in $\bf i$ (from the crossings at either
end of $c$) which is left unchanged by $B$ and thus determines a chamber
$(c',{\bf i'})$ in $\bf i'$. In this situation we also set
$(c,{\bf i})\sim (c',{\bf i'})$.

We extend this to a relation on chambers for arbitrary elements of
$\chi_n$ by specifying that $(c,{\bf i})\sim (c',{\bf i'})$
if there is a sequence $$(c,{\bf i})=(c_0,{\bf i}_0),(c_1,{\bf i}_1),\ldots
,(c_m,{\bf i}_m)=(c',{\bf i'})$$ with $(c_j,{\bf i}_j)\sim (c_{j+1},
{\bf i}_{j+1})$ for $j=0,1,\ldots ,m-1$. It is easy to see that $\sim$ is
an equivalence relation.
This relation has an alternative description:
\end{definitions}

\begin{lemma}
Suppose $(c,{\bf i})$ and $(c',{\bf i'})$ are chambers. Then
$(c,{\bf i})\sim (c',{\bf i'})$ if and only if
$l(c,{\bf i})=l(c',{\bf i'})$. In other words, two chambers are in the
same equivalence class if and only if they have the same chamber set.
\end{lemma}

\noindent {\bf Proof:} First, suppose that $(c,{\bf i})\sim (c',{\bf i'})$,
and that $\bf i$ and $\bf i'$ are related by a single braid relation.
By the definition of $\sim$ and the structure of the chamber ansatz, it
is clear that $l(c,{\bf i})=l(c',{\bf i'})$. It follows that this is
true for arbitrary $\bf i$ and $\bf i'$. The reverse implication
is~\cite[Lemma 2.7.2]{bfz1}.~$\square$

\section{Partial Quivers and Spanning Vectors}
In order to describe the spanning vectors of the Lusztig cones (that is, the
columns of the inverses of their defining matrices, $P_{\mathbf{i}}$),
we use the quiver-compatible reduced expressions. These reduced expressions
are spread sufficiently evenly throughout the set of all reduced expressions
that we can use them to obtain properties for arbitrary reduced expressions.

This approach leads to an alternative way of labelling the chambers in a
chamber ansatz, by `partial quivers', which we shall now describe.
We shall then use this together with a detailed analysis of the chamber
ansatz of a quiver-compatible reduced expression in order to give an explicit
description of the spanning vectors of an arbitrary Lusztig cone.

Thanks are due to R. W. Carter, who realised that it is often natural
to label a spanning vector (or chamber) by the set of quiver-compatible
reduced expressions for the longest word for which it appears; this led to
the idea of a partial quiver given below.

A {\em partial quiver} $P$ of type $A_n$ is
a copy of the Dynkin graph of type $A_n$ with some of the edges directed (at
least one, and possibly all), in such a way that the subgraph consisting of
the directed edges and vertices incident with them is connected.
We number the edges of a partial quiver from $2$ to $n$, starting at
the right hand end.
We write $P$ as a sequence of $n-1$ symbols, $L$, $R$ or $-$,
where $L$ denotes a leftward arrow, $R$ a rightward arrow, and $-$ an
undirected edge. Thus any partial quiver will be of form $---***---$, where
the $*$'s denote $L$'s or $R$'s.

If $P,P'$ are partial quivers we write $P'\geq P$ (or $P\leq P'$)
if every edge which is directed in $P$ is directed in $P'$ and is oriented in the
same way, and say that $P$ is a sub partial quiver of $P'$.
For example, $---LRLL-\ \leq RLRLRLLL$.

It is known that if $\bf i$ is compatible with a quiver $Q$
(in the sense of~\cite[\S4.7]{lusztig2}) then the
set of reduced expressions for $w_0$ compatible with $Q$ is precisely
the commutation class of $\chi_n$ containing $\bf i$.
We say that this commutation class is compatible with $Q$.
Berenstein, Fomin and Zelevinsky (see~\cite[\S4.4]{bfz1})
describe a method for constructing $\bf i\in \chi_n$ compatible with any given
quiver $Q$, as follows.

\begin{figure}
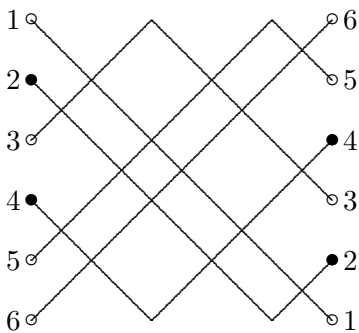

\begin{center}
\beginpicture
\setcoordinatesystem units <0.8cm,0.8cm>             
\setplotarea x from -7 to 7, y from 0 to 7       

\linethickness=0.5pt           

\multiput{$\circ$} at 1 1 *5 0 1 /
\multiput{$\circ$} at 6 1 *5 0 1 /

\put{$\bullet$} at 1 3 %
\put{$\bullet$} at 1 5 %
\put{$\bullet$} at 6 2 %
\put{$\bullet$} at 6 4 %

\put{$1$}[c] at 0.7 6 %
\put{$2$}[c] at 0.7 5 %
\put{$3$}[c] at 0.7 4 %
\put{$4$}[c] at 0.7 3 %
\put{$5$}[c] at 0.7 2 %
\put{$6$}[c] at 0.7 1 %

\put{$6$}[c] at 6.3 6 %
\put{$5$}[c] at 6.3 5 %
\put{$4$}[c] at 6.3 4 %
\put{$3$}[c] at 6.3 3 %
\put{$2$}[c] at 6.3 2 %
\put{$1$}[c] at 6.3 1 %

\setlinear \plot 1 6    6 1 / %

\setlinear \plot 1 5    5 1 / %
\setlinear \plot 5 1    6 2 / %

\setlinear \plot 1 4    3 6 / %
\setlinear \plot 3 6    6 3 / %

\setlinear \plot 1 3    3 1 / %
\setlinear \plot 3 1    6 4 / %

\setlinear \plot 1 2    5 6 / %
\setlinear \plot 5 6    6 5 / %

\setlinear \plot 1 1    6 6 / %

\endpicture
\caption{The arrangement Arr$\{2,4\}$.\label{arrangement}}
\end{center}
\end{figure}

Suppose $Q$ is a quiver of type $A_n$. Let $\Lambda\subseteq [2,n]$ be the
set of all edges of $Q$ pointing to the left. Berenstein, Fomin and
Zelevinsky construct an arrangement Arr$(\Lambda)$.
Consider a square in the plane, with horizontal
and vertical sides.
Put $n+1$ points onto the left-hand edge of the square, equally spaced,
numbered $1$ to $n+1$ from top to bottom, and
do the same for the right-hand edge, but number the points
from bottom to top. ${\rm Line}_h$ joins point $h$ on the left with
point $h$ on the right. For $h=1,n+1$, ${\rm Line}_h$ is a diagonal
of the square. For $h\in [2,n]$, ${\rm Line}_h$ is a union of two
line segments of slopes $1$ and $-1$. There are precisely
two possibilities for $\mbox{Line}_h$. If $h\in \Lambda$, the
left segment has slope $-1$, while the right one has slope $1$;
for $h\in [2,n]\setminus \Lambda$, it goes the other way round.
Berenstein, Fomin and Zelevinsky give the example of the arrangement 
for $n=5$ and $\Lambda=\{2,4\}$, which we show in Figure~\ref{arrangement} (points
corresponding to elements of $\Lambda$ are indicated by filled circles).
They show that $\bf i\in \chi_n$ is compatible with the quiver $Q$
if and only if the chamber ansatz CD($\bf i$) is isotopic to Arr$(\Lambda)$.

They also calculate the chamber sets of such a reduced expression
\cite[\S4.4.5]{bfz1}.

\begin{prop} \label{quiverchambers} (Berenstein, Fomin and Zelevinsky.) \\
Let ${\bf i}$ be a reduced expression compatible with a quiver $Q$.
Suppose that $(c,{\bf i})$ is a
chamber, and that strings $i$ and $j$ ($i<j$) cross immediately above it
in CD$(\bf i)$.
Then $l(ij):=l(c,{\bf i})$ is given in the following way:
Write out $Q$ and label the edges from $2$ to $n$, starting at the
right hand end. Define $X$ to be the subset of $[2,n]$ corresponding to
the $L$'s in $Q$. Then $l(ij)=Y_1\cup Y_2\cup Y_3$, where
$Y_1=[i+1,j-1]\cap X$,
$Y_2=[1,i-1]$ if $i\in [2,n]$ and $i\not\in X$, and is otherwise empty,
and $Y_3=[j+1,n+1]$ if $j\in [2,n]$ and $j\not \in X$, and is otherwise
empty.$~\square$
\end{prop}

Any given chamber set will arise in the chamber ansatz of a number of
reduced expressions. We now give a description of the quiver-compatible
reduced expressions whose chamber ansatz includes a fixed chamber set.
Suppose $(c,{\bf i})$ is a chamber. We denote by $Q(c,{\bf i})$ the
set of quivers $Q$ for which $l(c,{\bf i})$ occurs as the chamber
set of a chamber $(c',\bf i')$, where $\bf i'$ is compatible with $Q$.
If $P$ is a partial quiver, we define a subset $l(P)$ of $[1,n+1]$
in the following way.
Let $l_1(P)=\{j\in [2,n]\,:\,\mbox{edge\ $j$\ of\ $P$\ is\ an\ 
$L$}\}$.
If the rightmost directed edge of $P$ is an $R$, and this is
in position $i$, then let $l_2(P)=[1,i-1]$, otherwise the empty set.
If the leftmost directed edge of $P$ is an $R$, and this is
in position $j$, then let $l_3(P)=[j+1,n+1]$, otherwise the empty set.
Then put $l(P)=l_1(P)\cup l_2(P)\cup l_3(P)$. Note that for any partial
quiver $P$, $l(P)$ is always a chamber set, that is, it is a subset of
$[1,n+1]$ not of the form $[1,j]$ or $[j,n+1]$ for some $1\leq j\leq n+1$.
Thus $l$ is a map from chamber sets to partial quivers.

\begin{lemma}
The map $l$ is a bijection.
\end{lemma}
\noindent {\bf Proof:} First of all, note that the number of partial quivers
is the same as the number of chamber sets, i.e., $2^{n+1}-2(n+1)$, so it
is enough to show that $l$ is surjective. Let $S$ be a chamber set.
If $1\in S$, let $a$ be the minimal element of $[1,n+1]\setminus S$;
otherwise let $a$ be the minimal element of $S$. If $n+1\in S$,
let $b$ be the maximal element of $[1,n+1]\setminus S$; otherwise let $b$
be the maximal element of $S$. Let $P$ be a partial quiver with the edges $j$,
$a\leq j\leq b$ directed and the others undirected. For $a\leq j\leq b$,
let edge $j$ of $P$ be an $L$ if $j\in S$ and let it be an $R$ if not.
Then it can be seen from the definition of the map $l$ above that $l(P)=S$
and we are done.~$\square$

\begin{prop} \label{partiallabel}
Let $(c,{\bf i})$ be a chamber. Then there is a unique
partial quiver $P=P(c,{\bf i})$ such
that $Q(c,{\bf i})=\{Q\,:\,Q\mbox{\ is\ a\ quiver\ and\ }Q\geq P\}.$ 
Furthermore, if $(c,{\bf i})$ and
$(c',{\bf i}')$ are two chambers, then $l(c,{\bf i})=l(c',{\bf i'})$ if
and only if $P(c,{\bf i})=P(c',{\bf i}')$.
That is to say, their chamber sets are equal if and only if their partial
quivers are equal. We have
$l(P(c,{\bf i}))=l(c,{\bf i})$, and that every
partial quiver arises as $P(c,{\bf i})$ for some chamber $(c,{\bf i})$.
Thus the equivalence classes of
chambers are in one-to-one
correspondence with the partial quivers.
\end{prop}

\noindent {\bf Proof:}
We firstly claim that if we take a quiver $Q$, a reduced expression ${\bf i}$
compatible with $Q$ and a partial quiver
$P\leq Q$ then $l(P)=l(c,{\bf i})$ for some chamber $(c,{\bf i})$.
Fix a quiver $Q$ and let $P$ be a partial quiver with $P\leq Q$.
Suppose $P$ is of the form $----Y******X----$, where $X$, $Y$ and the
$*$'s are $L$'s or $R$'s and the $X$ appears in position $a$ and the
$Y$ appears in position $b$ ($a<b$), with the usual numbering.
If $X=L$ let $i$ be the position of the first $L$ in $Q$ strictly
to the right of $X$ (or $1$ if this doesn't exist).
If $X=R$ let $i=a$.
If $Y=L$ let $j$ be the position of the first $L$ in $Q$ strictly
to the left of $Y$ (or $n+1$ if this doesn't exist).
If $Y=R$ let $j=b$.
It is easy to check using Proposition~\ref{quiverchambers} that
$l(ij)=l(P)$; in fact the three components of $l(ij)$ and $l(P)$ coincide.

This deals with all of the cases except when
$P$ is of the form  $---R---$, where edge $i$ is an $R$.
Here $l(P)=[1,i-1]\cup [i+1,n+1]$.
But note that in the chamber ansatz for $Q$ string $i$ (corresponding to
an R) goes first up to the top of the diagram and then downwards
(by the description of Berenstein, Fomin and Zelevinsky). So at some
point it bounces off the top of the diagram. Therefore the chamber where it 
bounces will have chamber set precisely $[1,i-1]\cup [i+1,n+1]$. Thus the
claim is shown.

Since chambers of a reduced expression $\bf i$ correspond
to minimal pairs of equal letters, the chamber ansatz for $\bf i$ must have
precisely $\frac{1}{2}n(n-1)$ chambers. This is also the number of
partial quivers $P\leq Q$. Thus the chamber sets of $\bf i$ are precisely
the $l(P)$ for $P\leq Q$, $P$ a partial quiver (since $l$ is bijective).

Next, suppose that $(c,{\bf i})$ is a chamber, with chamber set
$l=l(c,{\bf i})$. Let $P$ be the partial quiver
such that $l(P)=l$. Then $l$ is a chamber set
of a reduced expression compatible with a quiver $Q$ if and only if
$P\leq Q$. Thus
$Q(c,{\bf i})=\{Q\,:\,Q\mbox{\ is\ a\ quiver\ and\ }Q\geq P\},$ as
required. It is clear that $P$ is unique with this property.
Also, we see that $l(P(c,{\bf i}))=l(c,{\bf i})$.
Furthermore, if $l(c,{\bf i})=l(c',{\bf i'})$, then since $P(c,{\bf i})$
is defined in terms of $l(c,{\bf i})$ it is clear that $P(c,{\bf i})=
P(c',{\bf i'})$. Conversely, if $P(c,{\bf i})=P(c',{\bf i'})=P$ then
applying $l$ to both sides gives $l(c,{\bf i})=l(c',{\bf i'})$, and we
are done.~$\square$

\begin{cor}
Let $Q$ be a quiver, and $\bf i$ a reduced expression compatible with
$Q$. Then the labels $P(c,{\bf i})$, for $(c,{\bf i})$ a chamber for $\bf i$,
are precisely the partial quivers $P\leq Q$.
\end{cor}

We next need to understand which crossings surround a chamber labelled by
a given partial quiver. We start with the following Lemma, which is easily
checked:


\begin{lemma} \label{chamberroots}
Suppose $X=(c,{\bf i})$ is a chamber in CD$({\bf i})$. Then let $\alpha$ and
$\beta$ be the positive roots labelling the crossings at the left and right
ends of $X$, and $\gamma_1,\gamma_2,\ldots ,\gamma_t$ the positive roots
labelling the crossings above and below $X$. Then
$\alpha+\beta=\sum_{i=1}^t \gamma_i$.~$\square$
\end{lemma}

In order to describe these crossings, we shall use the following definition.

\begin{definition}\rm
Let $2\leq a\leq b\leq n$, and let $Q$ be a quiver. We define the
$(a,b)$-sub partial quiver of $Q$ to be the partial quiver $P$ with
$P\leq Q$, leftmost directed edge numbered $b$ and rightmost directed edge
numbered $a$.
\end{definition}

\begin{lemma} \label{chambercrossings}
Let $Q$ be a quiver and ${\bf i}\in \chi_n$ a reduced expression compatible
with $Q$. Let $X=(c,{\bf i})$ be a chamber in CD$({\bf i})$, with
$P=P(c,{\bf i})$. Let $a$ be the position of the rightmost directed edge of
$P$, and $b$ the position of the leftmost directed edge.
We divide the possibilities for $P$ into four cases,
depending on the orientation of these edges of $P$:

%
%

$(1)\ \ \ \begin{array}{c} q \\ L \\ \ \end{array}\!\!\! RR\ldots R
\!\!\!\begin{array}{c} r \\ L \\ b\end{array}\ldots
\begin{array}{c} s \\ L \\ a \end{array}\!\!\! RR\ldots R
\!\!\!\begin{array}{c} p \\ L \\ \ \end{array}$

$(2)\ \ \ \begin{array}{c} q \\ L \\ \ \end{array}\!\!\! RR\ldots R
\!\!\!\begin{array}{c} r \\ L \\ b\end{array}\ldots
\begin{array}{c} p \\ R \\ a \end{array}\!\!\! LL\ldots L
\!\!\!\begin{array}{c} s \\ R \\ \ \end{array}$

$(3)\ \ \ \begin{array}{c} r \\ R \\ \ \end{array}\!\!\! LL\ldots L
\!\!\!\begin{array}{c} q \\ R \\ b \end{array}\ldots
\begin{array}{c} s \\ L \\ a \end{array}\!\!\! RR\ldots R
\!\!\!\begin{array}{c} p \\ L \\ \ \end{array}$

$(4)\ \ \ \begin{array}{c} r \\ R \\ \ \end{array}\!\!\! LL\ldots L
\!\!\!\begin{array}{c} q \\ R \\ b \end{array}\ldots
\begin{array}{c} p \\ R \\ a \end{array}\!\!\! LL\ldots L
\!\!\!\begin{array}{c} s \\ R \\ \ \end{array}$

In each case we define $4$ numbers, $p,q,r,s$, which number edges
of $Q$, defined by the above diagrams. For example, in the first case,
the leftmost and rightmost edges of $P$ are both $L$'s. In this case,
we define $r=b$, $s=a$, $q$ to be the position of the first $L$ in $Q$
appearing strictly to the left of the leftmost edge in $P$ and $p$ to be
the position of the first $L$ in $Q$ appearing strictly to the right of the
rightmost edge in $P$. The other cases are similar. We use the convention
that $Q$ has an edge in position $1$ with the same orientation as its
edge in position $a$, and that it has an edge in position $n+1$ with the
same orientation as its edge in position $b$ (thus $p,q,r,s$ are always
well-defined). Write $p(c,{\bf i})=p,$ $q(c,{\bf i})=q,$ $r(c,{\bf i})=r$ and
$s(c,{\bf i})=s$.

Make the convention that a crossing involving two equal strings means no
crossing. Then strings $p$ and $q$ cross immediately above $X$,
strings $q$ and $s$ cross to the left of $X$, strings $p$ and $r$ cross
to the right of $X$ and strings $r$ and $s$ cross below $X$.
\end{lemma}

\noindent {\bf Proof:} First of all, note that
each chamber in CD$({\bf i})$ has at
most one crossing bounding it on each side, left, right, above or below.
If $p\not=q$, the result is an easy consequence of
Proposition~\ref{quiverchambers} and the construction of Berenstein, Fomin and
Zelevinsky of quiver compatible reduced expressions.
If not, the result is easy to see, as the partial
quiver is a single $R$ in the $a$th position, so
$\ell(P)=[1,n+1]\setminus\{a\}$.~$\square$

Now we shall see how partial quivers can be used to describe the spanning
vectors of a Lusztig cone. Firstly, let us set up some notation.
Fix ${\bf i}\in \chi_n$. Then we know that
$$C_{\bf i}=\{{\bf a}\in \mathbb{Z}^k\,:\,P_{\bf i}{\bf a}\geq 0\},$$
where $P_{\bf i}$ is a $k\times k$ matrix as defined after
Lemma~\ref{addroots}, with inverse $Q_{\bf i}$ (over $\mathbb{Z}$). 
Note that $n(n-1)/2$
of the rows of $P_{\bf i}$ correspond to inequalities arising from minimal
pairs in $\bf i$. Each such minimal pair corresponds naturally to a chamber
$(c,{\bf i})$. Thus, given a chamber $(c,{\bf i})$,
there is a corresponding row of $P_{\bf i}$ and
therefore a corresponding column of $Q_{\bf i}$, that is,
a spanning vector of $C_{\bf i}$. We denote this spanning vector by
$v(c,{\bf i})$.
Similarly, the other $n$ rows of $P_{\bf i}$ correspond to inequalities
of the form $a_{\alpha_j}\geq 0$ for $\alpha_j$ a simple root, $1\leq j\leq
n$. We denote the corresponding spanning vectors by $v(j,{\bf i})$.

We shall usually regard elements $v$ of $C_{\mathbf{i}}$ as being indexed
by the positive roots (and thus by pairs of integers $p,q$ with
$1\leq p<q\leq n+1$). We shall occasionally need to consider such vectors
in the usual way, as elements of ${\mathbb{N}}^k$; we shall indicate this by
writing $v_{\mathbf{i}}$.

Suppose that $P$ is a quiver. We define a {\em sub partial quiver} $Y$ of a
quiver $P$ to be a {\em component} of $P$ if all of its edges are oriented
the same way and $Y$ is maximal in length with this property.
We say $Y$ has {\em type $L$} if its edges are oriented to the left, and
{\em type $R$} if its edges are oriented to the right.
For each component $Y$ of $P$, let $a(Y)$ be the position of the rightmost
directed edge of $Y$, and let $b(Y)$ be the position
of the leftmost directed edge of $Y$.
Let $v(Y)$ be the vector defined by setting, for $1\leq p<q\leq n+1$,
$v(Y)_{pq}=1$ if $1\leq p<a(Y)\leq b(Y)< q\leq n+1$ and $v(Y)_{pq}=0$,
otherwise.

\begin{theorem} \label{maintheorem} \ \\
(a) Let $(c,{\bf i})$ be a chamber. The spanning vector $v(c,{\bf i})$
depends only upon $P(c,{\bf i})$ (equivalently, upon the chamber set,
$\ell(c,{\bf i})$). For a partial quiver $P$ we choose a chamber
$(c,{\bf i})$ such that $P(c,{\bf i})=P$ and write $v(P)=v(c,{\bf i})$. \\
(b) For $j=1,2,\ldots ,n$, the spanning vector $v(j,{\bf i})$
depends only upon $j$; we denote it $v(j)$. \\
(c) For a partial quiver $P$, let
$$w(P)=\sum_{\mbox{\scriptsize $Y$\ a\ component\ of\ $P$}}v(Y).$$
Then, for each $1\leq p<q\leq n+1$, $v(P)_{pq}=\lceil
\frac{1}{2}w(P)_{pq} \rceil$ (recall that $\lceil\ \rceil$ rounds to the
nearest integer, rounding $\frac{1}{2}$ upwards).
Furthermore, for $j=1,2,\ldots n$, $v(j)_{pq}=1$ if
$1\leq p\leq j\leq j+1\leq q\leq n+1$ and $v(j)_{pq}=0$ otherwise.
\end{theorem}

{\bf Remark:}
Note that the theorem describes the columns of $P_{\mathbf{i}}^{-1}$, i.e.
the spanning vectors of $C_{\mathbf{i}}$, which is the problem we have set
out to solve.

\noindent {\bf Proof:}
We start by showing (a) and (b), i.e. that $v(c,\mathbf{i})$
(when its components are indexed by the positive roots) depends only on the
equivalence class of the chamber $(c,\mathbf{i})$ for the relation $\sim$
defined in~\ref{chamberequiv}, and that $v(j,\mathbf{i})$ (when its
components are indexed by the positive roots), depends only on $j$.
This reduces the question to quiver-compatible reduced words, which are
easier to deal with.

We start with two reduced expressions
${\bf i,\bf i'}\in\chi_n$, such that there is a short braid relation
$(i,j)\rightarrow (j,i)$ taking $\bf i$ to $\bf i'$.
Thus $\bf i'$ is obtained from
$\bf i$ by exchanging two of its entries. From the definition of the Lusztig
cone, it is clear that $P_{\bf i'}$ can be obtained from $P_{\bf i}$ by
exchanging the corresponding columns. Thus $Q_{\bf i'}$ can be obtained from
$Q_{\bf i}$ by exchanging the corresponding rows. This means that the
spanning vectors for $\bf i'$ are obtained from those for $\bf i$ by
exchanging the corresponding elements in the vectors, which means that
the spanning vectors for $\bf i'$, when labelled by the positive roots,
are the same as those for $\bf i$ (as the effect of the short braid relation
on the ordering of the positive roots defined by $\bf i$ is to exchange the
roots corresponding to the $i$ and the $j$).

Now suppose that $\bf i'$ is obtained from $\bf i$ by applying a long braid
relation $B:\,(i,j,i)\rightarrow (j,i,j)$. Since the pair of $i$'s in $\bf i$
is minimal, there is a corresponding chamber $X=(c,{\bf i})$ in
CD$({\bf i})$. Let $P'_{\bf i'}$ be the matrix obtained from $P_{\bf i'}$
by exchanging the columns corresponding to the pair of $j$'s in the long
braid relation in $\bf i'$.
Recall (see Definitions~\ref{chamberequiv}), that the effect of the long
braid relation on CD$({\bf i})$ is to map every chamber except
$X$ onto a chamber of CD$({\bf i'})$ with the same label; the chamber $X$
corresponds to a chamber $X'$ of CD$({\bf i'})$ with a label not equal
to any of the labels of the chambers of CD$({\bf i})$.
We thus have a correspondence between the chambers of CD$({\bf i})$ and the
chambers of CD$({\bf i'})$. If $(c,{\bf i})$ is a chamber in CD$({\bf i})$
we denote the corresponding chamber of CD$({\bf i'})$ by $(c,{\bf i})'$.
We denote the row of $P_{\bf i}$
corresponding to a chamber $Y$ of CD$({\bf i})$ by $r_Y$, and the row of
$P'_{\bf i'}$ corresponding to the corresponding chamber $Y'$ of
CD$({\bf i'})$ by $s_{Y'}$. We now prove the following claim.

{\bf Claim:} $P'_{\bf i'}$ can be obtained from
$P_{\bf i}$ by applying row operations which add or subtract $r_X$ to rows
other than $r_X$.

From the definition of the Lusztig cone we get that
$s_{X'}=r_X$. Let $Y\not= X$ be a chamber of CD$({\bf i})$ with corresponding
minimal pair of letters $p,\ldots ,p$ in $\bf i$. If this pair of letters
is to the left of the $(i,j,i)$ of $B$, or to the right of $B$,
then we have $s_{Y'}=r_Y$. The following can be checked easily.
Suppose first that the $(i,j,i)$ of $B$ lies
between the pair of $p$'s. If $|i-p|>1$ and $|j-p|>1$, then $s_{Y'}=r_Y$.
If $|i-p|=1$ and $|j-p|>1$ then $s_{Y'}=r_Y+r_X$. If $|i-p|>1$ and
$|j-p|=1$ then $s_{Y'}=r_Y-r_X$. Finally, if $p=i$ then $s_{Y'}=r_Y-r_X$,
and if $p=j$ then $s_{Y'}=r_Y+r_X$. Note also, that the rows corresponding to
simple roots are the same in $P_{\mathbf{i}}$ and in $P'_{\mathbf{i'}}$,
since the $j$ in $(i,j,i)$ can never correspond to a simple root in the
ordering defined by $\mathbf{i}$. Thus we see that the claim is proved.

It follows that there is a column operation on $Q_{\bf i}$ which changes
only the column corresponding to $X$ and which takes $Q_{\bf i}$ to the
inverse $Q'_{\bf i'}$ of $P'_{\bf i'}$. Now $Q'_{\bf i'}$ is the same
matrix as $Q_{\bf i'}$ except that the rows corresponding to the $i$'s in $B$
have been exchanged. Since the effect of $B$ on the ordering of the positive
roots defined by $\bf i$ is to exchange the roots corresponding to the $i$'s
in $B$, it follows that if $Y\not=X$ is a chamber of CD$({\bf i})$ then the
spanning vector corresponding to $Y$ is the same as the spanning vector
corresponding to $Y'$ (if both vectors are labelled by the positive roots),
and also, if $1\leq j\leq n$, then the spanning vector
corresponding to $(j,{\bf i})$ is the same as that for $(j,{\bf i'})$.

Thus statements (a) and (b) are true.
Note that every partial quiver $P$ satisfies $P\leq Q$ for some
quiver $Q$. So it is now enough to calculate the spanning vectors
for a reduced word $\bf i$ compatible with a quiver $Q$. We fix such a $Q$
and $\bf i$.

Fix $1\leq j\leq n$, and let $v$ be the vector defined by $v_{pq}=1$ when
$1\leq p\leq j\leq j+1\leq q\leq n+1$ and $v_{pq}=0$ otherwise. We shall
show that $v(j)=v$.
Let $P_{\bf i}$ be the defining matrix of $C_{\bf i}$
and let $r_l$, $1\leq l\leq n$, be the rows of $P_{\bf i}$
corresponding to the inequalities $v_{\alpha_l}\geq 0$ (where
$\alpha_1,\alpha_2,\ldots ,\alpha_l$ are the simple roots).
For each chamber $X$ in CD$({\bf i})$, let $r_X$ be the corresponding row of
$P_{\mathbf{i}}$.

It is clear that the product of $v$ with $r_l$ is
$\delta_{jl}$ for $1\leq l\leq n$, because $v_{l,l+1}=\delta_{jl}$.
Now let $X$ be a chamber of CD$({\bf i})$. To calculate the product of
$v$ with $r_X$ we have to calculate
$\sum_{i=1}^t v_{\gamma_i} -v_{\alpha} -v_{\beta}$, where the $\gamma_i$ are
the positive roots labelling the crossings immediately above and below $X$,
and $\alpha$ and $\beta$ are the crossings at the right and left ends of $X$.
But we know from Lemma~\ref{chamberroots} that $\sum_{i=1}^t\gamma_i=\alpha+
\beta$ from which we get that the product of $r_X$ with $v$
is zero. We conclude that $v(j)=v$ as required (as the spanning
vectors of $C_{\bf i}$ are the columns of the inverse of $P_{\bf i}$).

Suppose now that $P$ is a partial quiver with $P\leq Q$.
We know there is a chamber $Y$ of
CD$({\bf i})$ such that $P(Y)=P$. We have to show $v(P)=v$ where
$v$ is the vector as described in the Theorem.
It is clear that $v_{\alpha}=0$ if
$\alpha$ is simple, so we see that the product of $r_l$ with $v$ is
zero for $1\leq l\leq n$. So let $X$ be any chamber of CD$({\bf i})$.
We need to calculate $f(X,Y):=
\sum_{i=1}^t v_{\gamma_i}-v_{\alpha}-v_{\beta}$
(with $\gamma_i$, $\alpha$ and $\beta$ as above),
and show this is $\delta_{XY}$. Let $P'=P(X)$ be the partial quiver label
of $X$.

We have already calculated which crossings surround $X$, in
Lemma~\ref{chambercrossings}.
Let $p'=p(X)$, $q'=q(X)$, $r'=r(X)$ and $s'=s(X)$.
We must calculate $f(X,Y)=v_{p'q'}+v_{s'r'}-v_{s'q'}-v_{p'r'}$.
Note that $v_{ij}=\lceil w_{ij}/2\rceil$ where $w_{ij}$ is the number
of components of $P$ entirely contained in the $(i+1,j-1)$-sub partial quiver
of $Q$.
Let $a$ be the position of the rightmost directed edge of $P$, and let $b$
be the position of the leftmost directed edge of $P$, so that
$P$ is the $(a,b)$-sub partial quiver of $Q$.
Simlarly, let $a'$ be the position of the rightmost directed edge of $P'$,
and let $b'$ be the position of the leftmost directed edge of $P'$, so that
$P'$ is the $(a',b')$-sub partial quiver of $Q$.

We consider the four cases of Lemma~\ref{chambercrossings} for the partial
quiver $P'$. We consider first case (1), where the leftmost and rightmost
edges of $P'$ are oriented to the left. The diagram in the Lemma indicates
how to work out the values of $p',q',r'$ and $s'$. Let us recall it here:

$(1)\ \ \ \begin{array}{c} q' \\ L \\ \ \end{array}\!\!\! RR\ldots R
\!\!\!\begin{array}{c} r' \\ L \\ b'\end{array}\ldots
\begin{array}{c} s' \\ L \\ a' \end{array}\!\!\! RR\ldots R
\!\!\!\begin{array}{c} p' \\ L \\ \ \end{array}$

The first and last sequence of $R$'s may be empty, the other groupings not.
(This may happen, for example, if the leftmost edge of $P'$ and the edge
immediately to the left of it are both oriented to the left).

We note that:

(a) If $a>a'$, then $w_{s'r'}=w_{p'r'}$ and $w_{s'q'}=w_{p'q'}$. \\
(b) If $b<b'$, then $w_{s'r'}=w_{s'q'}$ and $w_{p'r'}=w_{p'q'}$.

In these cases, $f(X,Y)=0$ follows immediately from the definition, so we
can assume that $a\leq a'$ and $b\geq b'$.
Let us now consider the case where both the first and last sequences of
$R$'s occur.

It can be seen from the picture and the definition of $w_{ij}$ that:

(c) If $a=a'$, then $w_{s'r'}+1=w_{p'r'}$ and $w_{s'q'}+1=w_{p'q'}$. \\
(d) If $a<a'$, then $w_{s'r'}+2=w_{p'r'}$ and $w_{s'q'}+2=w_{p'q'}$. \\
(e) If $b>b'$, then $w_{s'r'}+2=w_{s'q'}$ and $w_{p'r'}+2=w_{p'q'}$. \\
(f) If $b=b'$, then $w_{s'r'}+1=w_{s'q'}$ and $w_{p'r'}+1=w_{p'q'}$.

We consider all possibilities, in each case computing
$w_{p'q'}$, $w_{s'r'}$, $w_{s'q'}$ and $w_{p'r'}$ in terms of $w_{s'r'}$ and
computing $f(X,Y)=
\lceil w_{p'q'}/2 \rceil + \lceil w_{s'r'}/2 \rceil -
\lceil w_{s'q'}/2 \rceil - \lceil w_{p'r'}/2 \rceil$.
We use the fact that $w_{s'r'}$ is odd.

Case (I) $a=a'$ and $b=b'$:
We have, by (c) and (f) above, that $w_{p'q'}=w_{s'r'}+2$, $w_{s'q'}=w_{s'r'}+1$ and $w_{p'r'}=w_{s'r'}+1$. Thus
$f(X,Y)=\lceil (w_{s'r'}+2)/2 \rceil + \lceil w_{s'r'}/2 \rceil -
\lceil (w_{s'r'}+1)/2 \rceil - \lceil (w_{s'r'}+1)/2 \rceil=1$.

Case (II) $a<a'$ and $b=b'$:
We have, by (d) and (f) above, that $w_{p'q'}=w_{s'r'}+3$, $w_{s'q'}=w_{s'r'}+1$ and $w_{p'r'}=w_{s'r'}+2$. Thus
$f(X,Y)=\lceil (w_{s'r'}+3)/2 \rceil + \lceil w_{s'r'}/2 \rceil -
\lceil (w_{s'r'}+1)/2 \rceil - \lceil (w_{s'r'}+2)/2 \rceil=0$.

Case (III) $a=a'$ and $b>b'$:
We have, by (c) and (e) above, that $w_{p'q'}=w_{s'r'}+3$, $w_{s'q'}=w_{s'r'}+2$ and $w_{p'r'}=w_{s'r'}+1$. Thus
$f(X,Y)=\lceil (w_{s'r'}+3)/2 \rceil + \lceil w_{s'r'}/2 \rceil -
\lceil (w_{s'r'}+2)/2 \rceil - \lceil (w_{s'r'}+1)/2 \rceil=0$.

Case (IV) $a<a'$ and $b>b'$:
We have, by (d) and (e) above, that $w_{p'q'}=w_{s'r'}+4$, $w_{s'q'}=w_{s'r'}+2$ and $w_{p'r'}=w_{s'r'}+2$. Thus
$f(X,Y)=\lceil (w_{s'r'}+4)/2 \rceil + \lceil w_{s'r'}/2 \rceil -
\lceil (w_{s'r'}+2)/2 \rceil - \lceil (w_{s'r'}+2)/2 \rceil=0$.

We next consider the possibility that the leftmost and rightmost
sequences of $R's$ do not occur, i.e. that $q'=r'+1$ or $p'=s'-1$.
We note that:

(g) If $q'=r'+1$ and $b=b'$ then (f) above still holds. \\
(h) If $q'=r'+1$ and $b>b'$ then $w_{s'r'}=w_{s'q'}$ and $w_{p'r'}=w_{p'q'}$.
\\
(i) If $p'=s'-1$ and $a=a'$ then (c) above still holds. \\
(j) If $p'=s'-1$ and $a<a'$ then $w_{s'r'}=w_{p'r'}$ and $w_{s'q'}=w_{p'q'}$.

As before, it follows immediately from the definition that $f(X,Y)=0$ in
cases $(h)$ and $(j)$.

Case (V) $q'=r'+1$ and $p'\not=s'-1$.
The case $b>b'$ is dealt with by (h) above.
If $b=b'$, then, by (g), the arguments in cases (I)
and (II) can be used to show $f(X,Y)=\delta_{XY}$.

Case (VI) $q'\not=r'+1$ and $p'=s'-1$.
The case $a<a'$ is dealt with by (j) above.
If $a=a'$, then, by (i), the arguments in cases (I)
and (III) can be used to show $f(X,Y)=\delta_{XY}$.

Case (VII) $q'=r'+1$ and $p'=s'-1$.
If $a<a'$ or $b>b'$, then by (h) or (j) above, we can see that $f(X,Y)=
\delta_{XY}$. We are left with the case where
$a=a'$ and $b=b'$, in which, by (g) and (i), the argument in case (I) can be
used to show $f(X,Y)=\delta_{XY}$.

We have thus covered all possibilities where the leftmost and rightmost
edges of $P'$ are oriented to the left --- case (1) of
Lemma~\ref{chambercrossings}. Computation of $f(X,Y)$ in case (2) of
Lemma~\ref{chambercrossings} can be
obtained from that in case (1) by
exchanging the roles of $p'$ and $s'$. The arguments are exactly the same,
except that in this case, $w_{p'r'}$ (playing the role of $w_{s'r'}$) is
{\em even}; also, it is the first sequence of R's or the final sequence of
L's which may or may not occur (i.e. we may have $q'=r'+1$ or $s'=p'-1$).
The argument for case (1) shows that
$v_{s'q'}+v_{p'r'}-v_{p'q'}-v_{s'r'}=-\delta_{XY}$ (note that $p$ and $s$
are swapped). The minus sign arises because $w_{p'r'}$ is even.
It follows that $f(X,Y)=v_{p'q'}+v_{s'r'}-v_{s'q'}-v_{p'r'}=
\delta_{XY}$. Similarly, case (3) is obtained by swapping $q$ and $r$ in
case (1), and case (4) is obtained by swapping both $q$ and $r$ and $p$ and
$s$ (in this last case, the argument for case (1) shows directly
that $f(X,Y)=\delta_{XY}$; here $w_{p'q'}$, playing the role of $w_{s'r'}$,
is again odd).

Thus, in all cases, we have seen that $f(X,Y)=\delta_{XY}$. We see that
this holds for any pair of chambers $X$ and $Y$, and thus the Theorem
is proved.~$\Box$

\noindent {\bf Acknowledgements} \\
Thanks are due to Professor R. W. Carter of the University of Warwick,
whose help, ideas and discussions were invaluable.
This paper was written while the author was an
EPSRC research assistant of Professor K.\,A. Brown at the University of
Glasgow, Scotland. Part of the work was completed while the author was
attending the programme `Representation Theory of Algebraic Groups and
Related Finite Groups' at the Isaac Newton Institute, Cambridge, in 1997.
The author would like to thank the referee for comments on an earlier
version of this manuscript.

\newcommand{\noopsort}[1]{}\newcommand{\singleletter}[1]{#1}

\end{document}